\documentclass{article}
\usepackage{amssymb}
\usepackage{amsmath}
\usepackage[english]{babel}
\author{E.Yu. Lerner}
\title{Multiplicative function instead of logarithm \\(an elementary approach)}
\date{}
\begin{document}

\maketitle

\vspace{-5mm}

V.I.~Arnold has recently defined the complexity of finite
sequences of zeroes and ones in terms of periods and preperiods of
attractors of a dynamic system of the operator of finite
differentiation. Arnold has set up a hypothesis that the
sequence of the values of the logarithm is most complicated or
almost most complicated. In this paper we obtain the necessary and
sufficient conditions which make this sequence (supplemented with
zero) most complicated for a more wide
class of operators. We prove that a sequence of values of a
multiplicative function in a finite field is most complicated
or almost most complicated for any operator divisible by the
differentiation operator.

\medskip

{\large \bf 1. Main results.}\ Let $x$ be a column consisting of $n$
elements of a finite field $\text{$\mathbb{F}_q$}$:
$$
x=\left(\begin{array}{c} x_1\\
\vdots \\
x_n
\end{array}\right),\quad x_i\in \text{$\mathbb{F}_q$}.
$$
Let $A$ stand for any linear operator defined on $\mathbb{F}_q^n$ which is
permutational with the cyclic shift operator $\delta$:
$$A\delta=\delta A, \quad \text{ where }\quad\delta\left(\begin{array}{c} x_1\\
x_2\\
\vdots \\
x_n
\end{array}\right)=\left(\begin{array}{c} x_2\\
\vdots \\
x_n\\
x_1
\end{array}\right).
$$
We call the operators with the mentioned property
translation-invariant. Let an operator $A$
map a vector $x$ into a certain vector $x'$: $x'=Ax$; let it do the vector $x'$ into that
$x''$, etc. Evidently, the sequence $x,x',x'',\ldots$, finally, becomes circular.

In~\cite{arnold} V.I.~Arnold defines the components of the initial vector
$x$ with the help of various functions $f$: $x_i=f(i)$. If $x$
belongs to the attraction domain of a cycle of the largest period (for
the given map $A$), and the preperiod of the sequence
$\tilde x=x,x',x'',\ldots$ is maximum possible, then the function $f$
is said to be most complicated. But if $x$ belongs to the attraction domain of the cycle of the maximal period and the preperiod is less than the maximal one by the unity, then $f$ is said to be almost most complicated.

In \cite{arnold} V.I.\,Arnold studies the properties of the finite differentiation operator~$\Delta$. It has the form $\Delta=\delta-I$, where
$I$ is the identical operator. In addition, Arnold sets up
hypotheses on the complexity of several functions for the map~$\Delta$.
In this paper we prove these hypotheses in more accurate statements. Note that the proofs are based only on the
elementary algebraic constructions, sometimes they are rather intricate. See~\cite{FAN} for a more elegant proof of Theorem~2.

Let $p$ be the characteristic of the field \text{$\mathbb{F}_q$} ($p$
is a prime number and $q=p^d$). If $n$ is divisible by~$p$, then the polynomial
functions are defined by the formula $f(i)=h\,(i\!\!\mod p)$, where $h$ is a
polynomial, whose coefficients belong to the field~\text{$\mathbb{F}_q$}. If $n$ is not divisible by $p$, then the polynomials are identical
constants from the field~\text{$\mathbb{F}_q$}. All polynomials belong to the attraction domain of the null vector. They do not include the most complicated functions of the operator~$\Delta$, except for the case, when $n=p^m$. In this particular case all functions are polynomials, and the polynomials of the degree~$m$ represent the most complicated functions of the differentiation operator (\cite{arnold}).

According to~\cite{arnold}, in the case, when~$q=2$, $A=\Delta$,
the attracting tree to any point of the cycles of the operator~$A$
has the same structure. For example, if~$n$ is not divisible
by~$p$, then the operator~$\Delta$ maps $q-1$ vectors which do not
belong to the limit cycle, into the same vector from the
attractor. In this case the graph of the map~$A$ is a set of limit
cycles (in~\cite{arnold} a cycle of the length $m$ is denoted
by~$O_m$). Each point of these cycles is accessible through a
``bunch'' which has~$q-1$ edges (with $d=1$, i.e., $q=p$, this
graph is denoted by~$T_p$). See Tables~1,~2 for the decomposition
of the graph of the map~$\Delta$ onto the connected components for
$q=2$ and~$q=3$ with prime odd values of~$n$ less than~50.
See~\cite{karpen},~\cite{tabl} for more detailed tables. The
computational algorithm is described in Section~3 (see also
\cite{Garber}).
\begin{center}
\begin{tabular}{|r||c|r|}
\multicolumn{3}{c}{Table~1. The graph of the map $\Delta:
\mathbb{F}_2^n\to\mathbb{F}_2^n$}\\
\hline n &the number of components& components of the graph of the map $\Delta$ ($q=2$)\\
\hline \hline
 3&2&$(O_3*T_2)+(O_1*T_2)$\\
 5&2&$(O_{15}*T_2)+(O_1*T_2)$\\
 7&10&$9(O_7*T_2)+(O_1*T_2)$\\
 11&4&$3(O_{341}*T_2)+(O_1*T_2)$\\
 13&6&$5(O_{819}*T_2)+(O_1*T_2)$\\
 \hline
 17&260&$256(O_{255}*T_2)+3(O_{85}*T_2)+(O_1*T_2)$\\
 19&28&$27(O_{9709}*T_2)+(O_1*T_2)$\\
 23&2050&$2049(O_{2047}*T_2)+(O_1*T_2)$\\
 29&566&$565(O_{475107}*T_2)+(O_1*T_2)$\\
 31&34636834&$34636833(O_{31}*T_2)+(O_1*T_2)$\\
 \hline
 37&21256&$21255 (O_{3233097}*T_2)+(O_1*T_2)$\\
 41&26214476&$26214400 (O_{41943}*T_2)+ 75 (O_{13981}*T_2)+(O_1*T_2)$\\
 43&805355524&$805355523 (O_{5461}*T_2)+(O_1*T_2)$\\
 47&8388610&$8388609 (O_{8388607}*T_2)+(O_1*T_2)$\\
\hline
\end{tabular}
\end{center}

\begin{center}
\begin{tabular}{|r||c|r|}
\multicolumn{3}{c}{Table~2. The graph of the map $\Delta:
\mathbb{F}_3^n\to\mathbb{F}_3^n$}\\
\hline n &the number of components& components of the graph of the map $\Delta$ ($q=3$)\\
\hline \hline
 5&2&$(O_{80}*T_3)+(O_1*T_3)$\\
 7&3&$2(O_{364}*T_3)+(O_1*T_3)$\\
 11&246&$243(O_{242}*T_3)+2(O_{121}*T_3)+(O_1*T_3)$\\
 13&20469&$20412(O_{26}*T_3)+56(O_{13}*T_3)+(O_1*T_3)$\\
 \hline
 17&194&$193(O_{223040}*T_3)+(O_1*T_3)$\\
 19&519&$518(O_{747916}*T_3)+(O_1*T_3)$\\
 23&177150&$177147 (O_{177146}*T_3)+2(O_{88573}*T_3)+(O_1*T_3)$\\
 29&82466&$82465 (O_{277412144}*T_3)+(O_1*T_3)$\\
 31&231435&$231434 (O_{889632172}*T_3)+(O_1*T_3)$\\
 \hline
 37&103053853533 &$103053850074 (O_{1456468}*T_3)+\phantom{11111111111111}$\\
& & $+3458 (O_{112036}*T_3)+(O_1*T_3)$\\
 41& 1853302661441610 & $1853302661441604 (O_{6560}*T_3)+\phantom{11111111111111}$\\
 &  & $+5(O_{1312}*T_3)+(O_1*T_3)$\\
\hline
 43& 121632015 & $121632014 (O_{899590375372}*T_3)+(O_1*T_3)$\\
 47&94143178830    & $94143178827 (O_{94143178826}*T_3)+\phantom{111111111111111}$\\
  & & $+2(O_{47071589413}*T_3)+(O_1*T_3)$\\
\hline
\end{tabular}
\end{center}

V.I.~Arnold has studied the complexity of various functions for the
operator~$\Delta$ with $q=2$. In particular, for the case, when $n+1$
is an odd prime number (we denote it by ~$r$), he has
considered the arithmetic logarithm defined by the formula
\begin{equation}
\label{1}
 f(i)=\left\{
\begin{array}{ll} 0,&\text{if $i$ is a quadratic residue modulo $r$;}\\
1,& \text{if $i$ is a quadratic nonresidue modulo $r$.}
\end{array}
\right.
\end{equation}
One can write the conditions ``$i$ is a quadratic residue modulo $r$'' and ``$i$ is a quadratic nonresidue modulo $r$'' in terms of the
Legendre symbol \text{\large $\left(\frac{i}{r}\right)$} as
$\text{\large $\left(\frac{i}{r}\right)$}=1$ and $\text{\large
$\left(\frac{i}{r}\right)$}=-1$, respectively.

The arithmetic logarithm is the most complicated or almost most complicated function of the differentiation operator for all $n\le 13$
(\cite{arnold}). Unfortunately, the hypothesis on the complexity of the logarithm for large dimensions appeared to be false (\cite{Garber}). Nevertheless, with the help of the quadratic Gauss sums in Section~2 we prove the following proposition (hereinafter the symbol $a\bot b$
introduced in \cite{knut} means that the $\text{GCD }(a,b)=1$).

\medskip

{\large Theorem 1.}\ {\it Let $n=r$ be an odd prime number.
Then the function $f$ defined by formula~(\ref{1}) for $1\le i\le n-1$
and extended as follows:
\begin{equation}
\label{2}
 f(n)=f(0)=0
\end{equation}
is the most complicated function of any translation-invariant
operator $A$ on $\mathbb{F}_q^n$ if and only if for certain integer~$k$ \\
either $n=4k+3$ and $q\,\bot\,(k+1)$ and $q\bot\,(2k+1)$,\\
or $n=4k+1$ and $q\,\bot\,2k$.
 }

\medskip
Since zero is an essentially singular point of the logarithmic function, it is not quite natural to redefine the logarithm by formula~(\ref{2}). Therefore in this paper we seek for complicated functions for which
the extension by the formula $f(0)=0$ is rather usual. The Legendre symbol itself represents such a function.

{\large Definition 1.}\ {\it Let $n$ be a prime number. A function
$f$ which maps $\{1,\ldots,n-1\}$ into \text{$\mathbb{F}_q$} is called
multiplicative, if $f(ij\!\!\mod n)=f(i)f(j)$ for any $i$,
$j$ from the definition domain, and $f$ differs from the identical zero.
 }

\medskip

Let~$S_q^m$ be the set of $m$th roots from the unity in the field
\text{$\mathbb{F}_q$}. Similarly to the case of multiplicative
characters (see \cite[Section 8.1]{Rou}) one can easily prove the following propositions:

1) $f(1)=1$;

2) all values $f(i)$ belong to $S_q^{n-1}$;

3) $\sum_{i=1}^{n-1}f(i)=0$, except for the case, when $f\equiv 1$.

We treat a multiplicative function which identically equals one as trivial. Property~2) implies that with $S_q^{n-1}=\{1\}$
no nontrivial multiplicative functions exist. This is true for
$(q-1)\,\bot\,(n-1)$, for example, for $q=2$. In general, the number of
various multiplicative functions equals the cardinal number of the
set $S_q^{n-1}$, i.e., the $\text{GCD }(n-1,q-1)$.
Really, one can easily prove that a multiplicative function is uniquely defined by its value at the point~$a$, where $a$
is a generatrix of the multiplicative group modulo~$n$.

The trivial function represents the unity of the group of multiplicative
functions, where the group operation is defined as the
componentwise multiplication. Therefore, $f^{-1}$ is a function
defined by the formula $f^{-1}(i)=(f(i))^{-1}$, $i=1,\ldots,n$.
Evidently, $f^{-1}(i)=f(i^{-1}\!\!\mod n)$, in particular,
$f^{-1}(n-1)=f(n-1)$. Similarly to the case of multiplicative characters
(see~\cite{Rou}) we redefine nontrivial multiplicative functions
by the formula~$f(0)=f(n)=0$.

\medskip

{\large Theorem 2.}\ {\it Let $n$ be an odd prime number; let
$f$ be a nontrivial multiplicative function, whose values belong to
\text{$\mathbb{F}_q$}; in addition, $n\ne p$, where $p$ is the characteristic
of the field~\text{$\mathbb{F}_q$}. Then $f$ is either the most complicated or almost most complicated function for any
translation-invariant operator $A: \mathbb{F}_q^n\to
\mathbb{F}_q^n$ representable in the form $A=B\Delta$.
 }

\medskip
Tables~1,~2 demonstrate that functions which differ from polynomials are not necessarily most complicated or almost most complicated functions of the operator~$\Delta$.

In the proof of Theorem~2 we use
the isomorphism of the algebra of cyclic matrices and the algebra of
polynomials of a variable ~$y$ defined modulo $y^n-1$. Here the matrix of the operator~$B$ with the first column
$
\left(\begin{array}{c} b_1\\
\vdots \\
b_n
\end{array}\right),\quad b_i\in \text{$\mathbb{F}_q$}$ corresponds to the polynomial ${\cal
B}(y)=\sum_{i=1}^{n} b_i y^{i-1}$.

\medskip

{\large Theorem 2$'$.}\ {\it Let conditions of Theorem~2 be fulfilled.
The function $f$ is the most complicated function of the operator~$A$
representable in the form $A=B\Delta$ if and only if
\begin{equation}
\label{th2} \text{\rm GCD }({\cal B}(y),\sum_{i=0}^{n-1} y^i)\ne
1,
\end{equation}
where ${\cal B}(y)$ is the polynomial, corresponding to the matrix of the
operator~$B$.}

\medskip

In particular, for $A=I \Delta$ we have the $\text{GCD
}(1,\sum_{i=0}^{n-1} y^i)=1$. Consequently, any nontrivial
multiplicative function is only the almost most complicated
function of the finite differentiation operator.

\medskip

{\large\bf 2. Proof of Theorem~1.}\  Let the symbol $A$ stand both for an operator and its matrix. Evidently, the
matrix of a translation-invariant operator~$A$ is a
cyclic matrix (a circulant \cite{Sachkov}), i.e., its
$i$th column $A_i$ satisfies the relation $A_{i}=\delta A_{i+1}$,
$i=1,\ldots,n-1$. This is equivalent to the fact, that each
row of a cyclic matrix is obtained from a previous one by a unit cyclic
shift to the right. Using the vector~$x$, we define a cyclic
matrix~$X$, whose first column coincides with~$x$. Evidently, the
equality $x'=Ax$ is equivalent to that $X'=AX$.

\medskip

{\large Lemma 1.} {\it Let
\begin{equation}
\label{xeqf} x=
\left(\begin{array}{c} f(1)\\
\vdots \\
f(n)
\end{array}\right),
\end{equation}
$f(i)\in \text{$\mathbb{F}_q$}$. A function $f$ is most
complicated for any translation-invariant
operator~$A$ if and only if the corresponding matrix
$X$ is nondegenerate over the field~\text{$\mathbb{F}_q$}.}

\medskip

Let us first prove the sufficiency. Assume that a matrix~$X$ is
nondegenerate; let $A^{m}X=A^l X$ with certain integer nonnegative~$m$~and~$l$. Then $A^{m}=A^l$ and the
equality $A^m Z=A^l Z$ is true for any cyclic
matrix~$Z$. This means that the function~$f$ is most complicated
for the map~$A$.

Let us prove the necessity ab contrario. Let a
vector $x$ be defined by formula (\ref{xeqf}), where $f$ is the most
complicated function of any
translation-invariant operator~$A$, in addition,
the corresponding matrix~$X$ is degenerate. If $X=\mathbb{O}$,
where $\mathbb{O}$ is the null matrix, then for all operators~$A$ the
sequence $y,Ay,\ldots$ immediately becomes cyclic with
any~$y\in\mathbb{F}_q^n$, what is impossible. Consequently,
$0<\text{rank } X<n$. Therefore, one can find a cyclic matrix $A$
such that $AX=\mathbb{O}$,
\begin{equation}
\label{rank} 0<\text{rank } A<n.
\end{equation}
Here, evidently, $(BA)X=\mathbb{O}$ for any matrix~$B$.

The fact that $f$ is the most complicated function for the maps $A$ and
$BA$ implies that $A^2=A$ and $(BA)^2=BA$. From the latter equality,
taking into account the commutativity of the multiplication for cyclic matrices,
we obtain that $AB=A$ for any nondegenerate cyclic
matrix~$B$. Consequently, $q=2$, $n=1$, what
contradicts~(\ref{rank}). Lemma~1 is proved.

\medskip

{\large Remark 1.} Lemma~1 implies that a finite analog of the
delta-function in the form
$$
f_k(i)=\left\{
\begin{array}{ll} 0,&\text{if $i=k$,}\\
1,& \text{otherwise,}
\end{array}
\right.
$$
is the most complicated function for any
translation-invariant operator. The action of the operator $\Delta$ onto this function is described in detail in~\cite{Garber}.

A similar proposition is true for linear operators defined on
$\mathbb{F}_q^n$ which correspond to a linear recurrent
correlation (see \cite[Section 8.2]{Lidl}).

\medskip

In order to complete the proof of Theorem~1, let us find the determinant of the
cyclic matrix~$X$, whose first column has form
(\ref{xeqf}), where the function~$f$ is defined by
formulas~(\ref{1}),~(\ref{2}).

It is well known (see \cite{Sachkov}) that eigenvalues of a cyclic
complex matrix are defined in terms of the elements of the first column
$f(j)$, $j=1,\ldots,n$ by the formulas
\begin{equation} \label{lambda}
\lambda_m=\sum_{j=1}^n f(j)\zeta^{j m},\
m=1,\ldots,n,\quad\text{where $\zeta=\exp(2\pi \imath/n).$}
\end{equation}
The determinant of the matrix equals the product of these values.
Note that a generalization of this formula onto the case,
when a cyclic matrix consists of arbitrary elements of a
finite field is considered in paper~\cite{duke}. However, it is more convenient for us to use formula~(\ref{lambda}).

The sum $g_m=\sum_{j=1}^n \text{\large $\left(\frac{j}{n}\right)$}
 \zeta^{jm}$ is said to be the
quadratic Gaussian sum~\cite{Rou}. Gauss has proved that
$g_m=\text{\large $\left(\frac{m}{n}\right)$}
 g_1$, where
$$
g_1= \left\{
\begin{array}{ll} \sqrt{n},&\text{if $n=4k+1$,}\\
\imath\sqrt{n},& \text{if $n=4k+3$.}
\end{array}
\right.
$$
We have
$$
2\lambda_m=-g_m + \sum_{j=1}^{n-1}\zeta^{jm}= \left\{
\begin{array}{ll} n-1,&\text{if $m=n$,}\\
-1-\text{\large $\left(\frac{m}{n}\right)$} g_1& \text{otherwise.}
\end{array}
\right.
$$
$$
\prod_{m=1}^n (2\lambda_m)=\left( g_1^2-1\right)^{(n-1)/2}(n-1).
$$
$$
\prod_{m=1}^n \lambda_m= \left\{
\begin{array}{ll} k^{(n-1)/2}\,2k,&\text{if $n=4k+1$,}\\
(k+1)^{(n-1)/2}\,(2k+1),& \text{if $n=4k+3$.}
\end{array}
\right.
$$
Evidently, the assertion of Theorem~1 follows from the latter formula and
Lemma~1.

\medskip

{\large \bf 3.~Corollaries of the isomorphism of the algebra of cyclic
matrices and the algebra of polynomials.}\newline Let
$\underline{a}=(\underline{a}_0,\underline{a}_1,\ldots,\underline{a}_{n-1})$,
$\underline{b}=(\underline{b}_0,\underline{b}_1,\ldots,\underline{b}_{n-1})$
be the first rows of cyclic matrices~$A$ and~$B$. Let $\underline{c}$
be the first row of the matrix~$C$ which represents the product of
matrices~$A$ and~$B$ (i.e., $C=A B$). We have
\begin{equation}
\label{matrPoly}
\underline{c}_{i}=\sum\limits_{j,k:\,j+k\!\!\!\!\mod n=i}
\underline{a}_{j} \underline{b}_{k},\quad i=0,\ldots,n-1.
\end{equation}
Equating the degrees which coincide modulo~$n$, we obtain a similar correlation for the coefficients of the
product of polynomials. To this end, suffice it to calculate the residues of the division by $y^n-1$, because
\begin{equation}
\label{mody} y^m\!\!\mod (y^n-1)=y^{m\!\!\mod n}.
\end{equation}

{\large Lemma~2.}\ {\it Let us associate any cyclic
matrix~$A$, whose first column is $
\left(\begin{array}{c} a_1\\
\vdots \\
a_n
\end{array}\right),\quad a_i\in\text{$\mathbb{F}_q$}$, with the polynomial ${\cal
A}(y)=\sum_{i=1}^{n} a_i y^{i-1}$. Let us define the multiplication of such
polynomials ${\cal A}(y)$ and ${\cal B}(y)$ by the formula $ {\cal
A}(y)\times{\cal B}(y)={\cal A}(y){\cal B}(y)\!\!\mod (y^n-1)$;
we understand the summation of polynomials and multiplication by a scalar value from \text{$\mathbb{F}_q$} as standard operations. Then the considered correspondence represents the isomorphism of the algebra of cyclic
matrices and that of polynomials.}
\medskip

Let us reformulate the assumption of Lemma~2,
replacing the first column with the first row, i.e., considering the polynomial $\underline{\cal
A}(y)=\sum_{i=0}^{n-1}\underline{a}_i y^i$ instead of that ${\cal
A}(y)$. Then, evidently, the assertion of the lemma follows from
formula~(\ref{matrPoly}). However, it is more convenient to prove
theorems~2 and~2$'$, considering the columns.
The connection between the elements of the first row and the first column is defined by the formula $\underline{a}_i=a_{(n+1-i)\!\!\!\mod\! n}$. In order to prove the lemma, suffice it to note that the map ${\cal
A}(y)\to \underline{\cal A}(y)=\sum_{i=0}^{n-1}
a_{(n+1-i)\!\!\!\mod\! n}\ y^i=y^n {\cal A}(1/y)\!\!\mod(y^n-1)$
represents the automorphism of the mentioned algebra of polynomials.

\medskip
{\large Lemma~3.}\ {\it Let $F(y)=\sum_{i=0}^{n-1} f(i)y^i$,
$f(i)\in\text{$\mathbb{F}_q$}$, $n\bot q$. If
$F(y)\bot\sum_{i=0}^{n-1}y^i$, then the function $f$ is the most complicated or almost most complicated function of any operator~$A$
representable in the form
\begin{equation} \label{ABDelta} A=B\Delta. \end{equation} In addition, if
$\sum_{i=1}^n f(i)=0$, then the function~$f$ is the most complicated
function of the operator~$A$ if and only if
condition~(\ref{th2}) is fulfilled.}
\medskip

{\large Proof.}\ Let~$x$ have form~(\ref{xeqf}). Due to
Lemma~2 the equality $A^m x=A^l x$ is equivalent to the correlation ${\cal
A}^m(y) \sum_{i=1}^n f(i)y^{i-1}\!\!\mod (y^n-1)=$
${\cal A}^l(y) \sum_{i=1}^n f(i)y^{i-1}\!\!\mod (y^n-1)$ or
\begin{equation} {\cal A}^m(y) F(y)\!\!\mod (y^n-1)={\cal A}^l(y)
F(y)\!\!\mod (y^n-1), \label{AmFeqAlF}
\end{equation}
where ${\cal
A}(y)=\sum_{i=1}^{n} a_i y^{i-1}$, $
\left(\begin{array}{c} a_1\\
\vdots \\
a_n
\end{array}\right)
$ is the first column of the matrix~$A$. In view of the decomposition
$y^n-1=(y-1)\sum_{i=0}^{n-1}y^i$, equality~(\ref{AmFeqAlF}), in
turn, is equivalent to the system
\begin{equation}\label{system}
\left\{\begin{array}{c} {\cal A}^m(y)
F(y)\!\!\!\mod (y-1)={\cal A}^l(y) F(y)\!\!\!\mod (y-1),\\
{\cal A}^m(y)F(y)\!\!\!\mod \sum_{i=0}^{n-1}y^i={\cal A}^l(y)
F(y)\!\!\!\mod \sum_{i=0}^{n-1}y^i.
\end{array}\right.
\end{equation}
Here we use the fact that $n\bot q$, consequently,
$y=1$ is not a root of the polynomial $\sum_{i=0}^{n-1}y^i$
in~\text{$\mathbb{F}_q$}, i.e., $\sum_{i=0}^{n-1} y^i\ \bot\
(y-1)$.

Formula~(\ref{ABDelta}) means that ${\cal A}(y)$ admits the
representation ${\cal A}(y)={\cal B}(y)(y^{n-1}-1)\!\!\mod(y^n-1)={\cal
B}(y)y^{n-1}(1-y)\!\!\mod(y^n-1)$, whence ${\cal A}(y)$ is divisible
by~$(y-1)$. Therefore, the first correlation in~(\ref{system})
is an identity with all~$m,l\ge 1$.

Let us find the period and the preperiod of the sequence~$\tilde
x=x,x',x'',\ldots$ for the most complicated function~$f$. Let $f$ be the
delta-function mentioned in Remark~1. For it we have~$F(y)=1$.
System (\ref{system}) implies that the preperiod of the corresponding
sequence~$\tilde x$ is defined as $\max(1,m')$, where
$m'$ is the preperiod of the sequence
\begin{equation}\label{series} {\cal A}^m(y)\!\!\!\mod
\sum_{i=0}^{n-1} y^i,\quad m=1,2,\ldots.\end{equation}
Here the period of the sequence $\tilde x$ equals that of
sequence~(\ref{series}).

Let us now estimate the preperiod and the period of the sequence~$\tilde x$
for a function~$f$ such that $F(y)\ \bot\
\sum_{i=0}^{n-1}y^i$. Since $F(y)$ is invertible in the algebra of
polynomials modulo~$\sum_{i=0}^{n-1}y^i$, the second correlation
in~(\ref{system}) in this case is equivalent to the equality $ {\cal
A}^m(y)\!\!\mod \sum_{i=0}^{n-1}y^i={\cal A}^l(y) \!\!\mod
\sum_{i=0}^{n-1}y^i$. So the preperiod and the period of the
sequence~$\tilde x$ are not less than the preperiod and the
period of sequence~(\ref{series}). The first part of Lemma~3
is proved.

If $\sum_{i=1}^n f(i)=0$, then $F(y)$ is divisible by $(y-1)$. Then
the first correlation in (\ref{system}) is true with all $m,l$.
Consequently, in this case the function $f$ corresponds to the
sequence $\tilde x$, whose preperiod equals $m'$. The function $f$
is most complicated if and only if $m'\ge 1$. This is
equivalent to the fact that ${\cal A}(y)$ is noninvertible in the algebra modulo
$\sum_{i=0}^{n-1}y^i$, i.e., to condition~(\ref{th2}). Lemma~3 is proved.

\medskip

{\large Remark 2.} An application of Lemma~2 enables us to define the structure
and the number of connected components for any
translation-invariant operator~$A$. To this end, suffice it to
decompose $y^n-1$ onto irreducible polynomials in the field
\text{$\mathbb{F}_q$}:
\begin{equation}
\label{circle} y^n-1=\prod_{j=1}^m P_j^{\beta_j}(y).
\end{equation}

Let us consider the polynomial ${\cal A}(y)$ (hereinafter in Remark~2 one can consider the polynomial $\underline{\cal A}(y)$ instead of that ${\cal
A}(y)$; it does not affect the result). Let ${\cal A}(y)$ be coprime with the polynomials
$P_1(y),\ldots,P_k(y)$ and let it be divisible by the polynomials $P_{k+1}(y),\ldots
P_m(y)$. Then the structure of the attracting tree to any point of the
attractor is the same for all points. Namely, it represents a tree with~$l$
levels, where $l$ is the minimal natural number such that ${\cal
A}^l(y)$ is divisible by $\prod_{j=k+1}^m P_j^{\beta_j}(y)$. In addition,
all vertices of the tree located at the same level have the same
number of sons. One can find it, using the fact that the
total amount of vertices at the first~$i$ levels equals $q^{r_i}$, where
$r_i$ is the degree of the polynomial the $\text{GCD }({\cal A}^i(y),\prod_{j=k+1}^m
P_j^{\beta_j}(y))$.

In order to determine the number of cycles of various lengths, suffice it to find natural numbers $s_{ij}$ which represent the orders of the polynomial~${\cal A}(y)$ considered
as an element of a multiplicative group modulo $P_j^i(y),
j=1,\ldots,k,\quad i=1,\ldots,\beta_j$. This order has to be a
divisor of the total amount of elements in the group, i.e., the number
$d_{ij}=q^{i\ \text{\rm deg } P_j}-q^{(i-1) \text{\rm deg } P_j}$,
where $\text{deg } P$ is the degree of the polynomial~$P$. To put it more precisely, the order
$s_{ij}$
is the least of the divisors such that
$ {\cal A}^{s_{ij}}(y)\!\!\!\mod P_j^i(y)=1. $ The connected
components of the operator~${\cal A}$ are
\begin{equation}\label{structure}
\prod_{j=1}^k\left(O_1+\sum_{i=1}^{\beta_j}\frac{d_{ij}}{s_{ij}}O_{s_{ij}}\right)*T,
\end{equation}
where~$T$ is the tree described above, $O_m$ is the cycle of the length~$m$, and
the product of different cycles is defined by the formula
\begin{equation}\label{ruleGCDLCM} O_m O_l=\text{GCD }(m,l)\
O_{{\rm LCM }(m,l)}. \end{equation}

\medskip

{\large An example.} Let $q=3$, $n=12$, $A=\Delta$. In the field \text{$\mathbb{F}_3$}, \ $ y^{12}-1=(y+2)^3(y+1)^3(y^2+1)^3,\quad
\underline{\cal A}(y)=y-1\equiv y+2$. Therefore, the tree $T$
has 3 levels, it consists of 27 vertices, each vertex, except for
leaves and the root, has 3 sons; the root has 2 sons.
Following~\cite{arnold}, we denote this tree by $T_{27}$.

The numbers of elements in the multiplicative groups modulo
$(y+1)$, $(y+1)^2$, $(y+1)^3$ in the field \text{$\mathbb{F}_3$}\ are equal,
correspondingly, to $3-1=2,\quad 3^2-3=6,\quad 3^3-3^2=18$. Further,
$(y-1)\!\!\mod(y+1)=1,\quad (y-1)^2\!\!\mod(y+1)^2\ne 1,\quad
(y-1)^2\!\!\mod(y+1)^3\ne 1$, but
$(y-1)^3\!\!\mod(y+1)^2=(y-1)^3\!\!\mod(y+1)^3=1$. Therefore, the orders
$(y-1)$ in multiplicative groups modulo $(y+1),\quad
(y+1)^2,\quad (y+1)^3$ are equal, correspondingly, to $1,3,3$. Thus, the multiplier $(y^2+1)^3$ in product (\ref{structure})
corresponds to the sum $O_1+2 O_1+6/3\ O_3+18/3\ O_3=3\ O_1+8\ O_3. $

The numbers of the elements in multiplicative groups modulo
$(y^2+1)$, $(y^2+1)^2$, $(y^2+1)^3$ in the field \text{$\mathbb{F}_3$}\
are equal, correspondingly, to $9-1=8,\quad 9^2-9=72,\quad 9^3-9^2=648$.
All possible divisors of these numbers are
\begin{eqnarray}
& &1,2,4,8;\label{pow1}\\
& &1,2,3,\ldots,24,36,72;\label{pow2}\\
& &1,2,3,\ldots,24,\ldots,216,324,648.\label{pow3}
\end{eqnarray}
Raising $(y-1)$ to powers (\ref{pow1}), (\ref{pow2}), (\ref{pow3})
modulo $(y^2+1)$, $(y^2+1)^2$, $(y^2+1)^3$ until we obtain
the unit, we determine the orders of $(y-1)$ in the corresponding
multiplicative groups. They equal $8,\ 24$, and~$24$. Therefore, the multiplier $(y+1)^3$ in product (\ref{structure})
corresponds to the sum $O_1+O_3+3\ O_{24}+27\ O_{24}=O_1+O_3+30\ O_{24}$.

Note that one can calculate the order of $(y-1)$ in a group of
$m$ elements in a different way. Let $m=\prod p_i^{b_i}$, where $p_i$
are prime numbers. In order to find the order, suffice it to divide $m$ consecutively (evidently, no more than $b_i$ times) onto prime numbers~$p_i$ which produce this number, until $(y-1)$ in the proper degree equals the unit. For example, the order of $(y-1)$ in the
multiplicative group modulo $(y^2+1)^2$ equals 24, because
$72=2^3 3^3$ and $(y-1)^{72/2}\!\!\mod(y^2+1)^2\ne 1$,
$(y-1)^{72/3}\!\!\mod(y^2+1)^2=1$,
$(y-1)^{24/3}\!\!\mod(y^2+1)^2\ne 1$.

According to rule (\ref{ruleGCDLCM}),
$$(3\ O_1+8\ O_3)(O_1+O_8+30\ O_{24})=3\ O_1+8\ O_3+3\ O_8+(8+3\times 30+24\times 30)
O_{24}.
$$
Thus, the graph of the map $\Delta:
\mathbb{F}_3^{12}\to\mathbb{F}_3^{12}$ admits the following decomposition onto the connected components:
$$
3(O_1*T_{27})+8(O_3*T_{27})+3(O_8*T_{27})+818(O_{24}*T_{27}).
$$

See \cite{tabl} (see also http://kek.ksu.ru/kek2/myArnold.htm) for
the decomposition of the graph of the map  $\Delta:
\mathbb{F}_q^n\to\mathbb{F}_q^n$ for all $n\le 300$ with $q=2$ and for
$n\le 150$ with $q=3$, as well as the Wolfram Research Mathematica program which performs this decomposition.

\medskip

{\large Remark 3.} One can simplify the algorithm described above and implemented in~\cite{tabl}, using several facts mentioned
in~\cite{Garber}. One can easily prove (see proposition~12 in
paper~\cite{Garber}, cf. with \cite[theorem~2.42]{Lidl}) that all
degrees~$\beta_j$ in formula (\ref{circle}) are equal to $p^m$, where $n=p^m
n'$, $n'\bot n$. In addition, lemma~5 from paper \cite{Garber} (cf.
with \cite[theorem~3.8]{Lidl}) implies that the order of $(y-1)$ modulo
$P^j(y)$ equals the order modulo $P(y)$ multiplied by
$p^k$, where $k=\lceil\log_p j\,\rceil$.

\medskip

{\large \bf 4. Proof of theorems~2~and~2$'$.}\ In accordance with
Lemma~3 and property~3) of nontrivial multiplicative
functions, suffice it to prove that for any function $f$ with the mentioned property,
 \begin{equation} \label{mainforth2}
 \mbox{GCD }(\sum_{i=1}^{n-1}f(i)y^i,\sum_{i=0}^{n-1}y^i)=1.
 \end{equation}
Let $F(y)$ stand for the polynomial $\sum_{i=1}^{n-1}f(i)y^i$, let
$F^-(y)$ stand for $\sum_{i=1}^{n-1}f^{-1}(i)y^i$, and let $H(y)$ denote the
product $F(y)F^-(y)$. Put
\begin{equation}\label{defh} H(y)\!\!\mod (y^n-1)=\sum_{i=0}^{n-1} h(i)y^i.
\end{equation}
Let us prove that
\begin{equation}\label{hequal}
h(1)=h(2)=\ldots=h(n-1).
\end{equation}
For positive integer
$k$ let us define the operator $\cdot_k$ which maps the polynomial
$S(y)=\sum_{i=0}^m s(i)y^i$ onto that $S_k(y)=\sum_{i=0}^m
s(i)y^{ik}$. Using (\ref{mody}), we obtain that for any $k$ which is
coprime to $n$,
$$
(f(k)F_k(y))\!\!\mod(y^n-1)=\sum_{i=1}^{n-1}f(i)f(k)\,
y^{ik\!\!\mod n}=\sum_{i=1}^{n-1}f(ik)\, y^{ik\!\!\mod n}=F(y).
$$
Analogously, $(f^{-1}(k)F^-_k(y))\!\!\mod(y^n-1)=F^-(y)$. Hence
$$H_k(y)\!\!\!\!\mod(y^n-1)=F_k(y)F^-_k(y)\!\!\!\!\mod(y^n-1)=
$$
$$
=(f(k)F_k(y))(f^{-1}(k)F^-_k(y))\!\!\!\!\mod(y^n-1)=F(y)F^-(y)\!\!\!\!\mod(y^n-1)
=\sum_{i=0}^{n-1} h(i)y^i.
$$
On the other hand, due to property~(\ref{mody}),
$$
H_k(y)\!\!\mod(y^n-1)=\sum_{i=0}^{n-1} h(i)y^{ik\!\!\mod n}.
$$
Equality~(\ref{hequal}) is proved.

Let us now prove that
\begin{equation}\label{hnequal} h(0)-h(1)\ne 0. \end{equation}
Let us calculate $h(0)$, i.e., the coefficient at $y^n$ in $H(y)=F(y)F^-(y)$.
It equals
$$
\sum_{i=1}^{n-1} f(i)f^{-1}(n-i)=\sum_{i=1}^{n-1}
f(i)f^{-1}(i)f^{-1}(-1)=(n-1)f^{-1}(-1)=(n-1)f(-1)
$$
(hereinafter $f(-1)\equiv f(n-1)$). Since $F(y)$ is divisible by
$(y-1)$, we conclude that $(F(y)F^-(y))\!\!\mod(y^n-1)$ is also divisible by
$(y-1)$, i.e., $\sum_{i=0}^{n-1}h(i)=0$. Due to~(\ref{hequal}) we have
$\sum_{i=1}^{n-1}h(i)=h(1)(n-1)$, then
$\sum_{i=0}^{n-1}h(i)=(n-1)(f(-1)+h(1))$, whence $h(1)=-f(-1)$,
consequently, $h(0)-h(1)=n f(-1)$. Since $n\bot p$ and $f(-1)\ne
0$, relation~(\ref{hnequal}) is true.

Let us now find $(F(y)F^-(y))\!\!\mod\sum_{i=0}^{n-1}y^i$, i.e.,
$H(y)\!\!\mod\sum_{i=0}^{n-1}y^i$. Evidently,
\begin{equation}\label{polyfinish}
H(y)\!\!\mod\sum_{i=0}^{n-1}y^i=\sum_{i=0}^{n-1}h(i)
y^i\!\!\mod\sum_{i=0}^{n-1}y^i.
\end{equation}
Due to~(\ref{hequal}) polynomial (\ref{polyfinish}) is an identical constant, and in view of~(\ref{hnequal}) this constant differs from zero. This means
that $F(y)$ is an invertible element of the algebra of polynomials
modulo~$\sum_{i=0}^{n-1}y^i$, consequently, formula~(\ref{mainforth2}) is true. Theorems~2 and 2$'$ are proved.

\medskip

{\large Remark 4.}\ Theorem~2 is false if the operator~$A$ does not admit the representation $A=B\Delta$.
The operator $A=I+\delta+\delta^2+\ldots+\delta^{n-1}$ gives a disproving example for the mentioned incorrect statement of the theorem.
Any nontrivial multiplicative function in this case defines a
sequence, whose period is less than the maximal one.


\begin{thebibliography}{9}
\bibitem{Rou}
{\it K.~Ireland, M.~Rosen.}\ A Classical Introduction to Modern
Number Theory, Springer, 1990.
\bibitem{arnold}
{\it V.\,I.~Arnold.}\ Complexity of finite sequences of zeros and
ones and geometry of finite spaces of functions. Functional
Analysis and Other Mathematics, {\bf 1}, No.~1, 1--18, (2006).
\bibitem{knut}
{\it R.\,L.~Graham, D.\,E.~Knuth, O.~Patashnik.}\ Concrete
Mathematics, Reading, Massachusetts: Addison-Wesley, 1994.
\bibitem{Lidl}
{\it R.\,Lidl, H.\,Niederreiter.}\ Finite Fields, Vol. 20 in the
Encyclopedia of Mathematics and its Applications, Addison-Wesley,
1983.
\bibitem{Sachkov}
{\it V.\,N.~Sachkov, V.E.~Tarakanov.}\ Combinatorial calculus of
nonnegative matrices. TVP, Moscow, 2000.
\bibitem{Garber}
{\it A.\,I.~Garber.}\ Graph of difference operators for $p$-ary
sequences. Functional Analysis and Other Mathematics, {\bf 1},
No.~2, 179--195, (2006).
\bibitem{karpen}
{\it O.\,N.~Karpenkov.}\ Om examples of difference operators for
$\{ 0,1\}$-valued functions over finite sets. Functional Analysis
and Other Mathematics, {\bf 1}, No.~2, 197--202, (2006).
\bibitem{tabl} {\it
E.\,Yu.~Lerner.}\ Tables of graphs of binary and ternary
sequence's dif\-fe\-ren\-ti\-ation. Preprint.
http://arxiv.org/abs/0704.2947.
\bibitem{FAN}
{\it E.\,Yu.~Lerner.}\ Complexity of prime finite differences. To
be appear.
\bibitem{duke}
{\it O.~Ore}\ Some studies on cyclic determinant. Duke
Mathematical Journal, {\bf 18}, No.~2, 343--354, (1951).
\end{thebibliography}
\end{document}